\documentclass{article} %
\usepackage{iclr2022_conference,times}
\PassOptionsToPackage{numbers, compress}{natbib}

\usepackage{amsmath, mathtools}
\usepackage{amsfonts}
\usepackage{amssymb}
\usepackage{graphicx}
\usepackage{verbatim}
\usepackage{mathpazo}
\usepackage{enumerate}
\usepackage{xcolor}
\usepackage{fourier}
\usepackage{physics}
\usepackage{amsmath,amssymb,latexsym}
\usepackage{subfigure}
\usepackage{graphicx}
\usepackage{booktabs}
\usepackage{epstopdf}
\usepackage{subfigure}

\usepackage{threeparttable}

\usepackage[framemethod=TikZ]{mdframed}
\usepackage{tcolorbox}

\usepackage{dsfont}
\usepackage[english]{babel}
\usepackage{mathtools}

\usepackage{url}
\usepackage{hyperref}

\DeclarePairedDelimiterX{\inp}[2]{\langle}{\rangle}{#1, #2}

\newtheorem{remark}{Remark}
\newtheorem{assumption}{Assumption}

\usepackage{xcolor}

\usepackage{graphicx}          

\def\p{\boldsymbol{p}}
\def\q{\boldsymbol{q}}

\def\x{\boldsymbol{x}}

\def\f{\boldsymbol{f}}
\def\F{\boldsymbol{F}}

\def\r{\boldsymbol{r}}

\def\u{\boldsymbol{u}}

\def\1{\boldsymbol{1}}

\newtheorem{problem}{Problem}
\newcounter{example}[section]
\newenvironment{example}[1][]{\refstepcounter{example}}{\medskip}

\title{Parametrization Approach for Real-Time Generation of Minimum-Effort Trajectories via  Neural Network}

\author{Han Wang \& Zheng Chen \\
	School of Aeronautics and Astronautics, Zhejiang University\\
	\texttt{z-chen@zju.edu.cn}
}

\iclrfinalcopy %

\begin{document}
\maketitle
\begin{abstract}                
This paper is concerned with real-time generation of optimal flight trajectories for Minimum-Effort Control Problems (MECPs), which is fundamentally important for autonomous flight of aerospace vehicles. Although existing optimal control methods, such as indirect methods and direct methods, can be amended to solve MECPs, it is very challenging to obtain, in real time, the solution trajectories since those methods suffer the issue of convergence. As the artificial neural network can generate its output within a constant time, it has been  alternative for real-time generation of optimal trajectories in the literature. The usual way is to train neural networks by solutions from indirect or direct methods, which, however, cannot ensure sufficient conditions for optimality to be met. As a result, the trained neural networks cannot be guaranteed to generate at least locally optimal trajectories. To address this issue, a parametrization approach is developed in the paper so that not only necessary but also sufficient conditions are embedded into a parameterized set of differential equations. This allows generating the dataset of at least locally optimal trajectories through solving some initial value problems. Once a neural network is trained by the dataset constructed by the parametrization approach, it can not only generate optimal trajectories within milliseconds but also ensure the generated trajectories to be at least locally optimal, as finally demonstrated by two conventional MECPs in aerospace engineering. 
\end{abstract}



\section{Introduction}\label{SE:Introduction}

It is widely required in aerospace engineering to generate flight trajectories in real time while minimizing control effort, such as  minimum-energy control of satellites \cite{athans1963time,horri2011energy}, trajectory optimization for hypersonic vehicles \cite{lu2006adaptive}, and  missile guidance \cite{chen2019nonlinear,Kun:2021}, just to name a few. Essentially, real-time generation of minimum-effort trajectories is related to solving some Minimum-Effort Control Problems (MECPs) within a constant time by onboard computers. As the MECP belongs to a specific class of optimal control problems, existing optimal control methods can be readily applied, and they can be generally categorized into two classes  \cite{doi:10.1137/1.9780898718577}: 1) indirect methods and 2) direct methods.

Indirect methods are based on first-order necessary conditions of optimality from calculus of variations or from Pontryagin's Maximum Principle (PMP). These conditions allow to transform an optimal control problem into a 2-Point Boundary Value Problem (2PBVP) or Multi-Point Boundary Value Problem (MPBVP), which can be solved by some rooting-finding or shooting-like methods. The optimal solution, obtained from indirect methods, is accurate. However, it is well known that root-finding or shooting-like methods are hard to converge unless a good guess in a small domain of initial costate is provided. In order to improve the convergence of indirect methods, continuation or homotopy methods have been introduced into shooting methods in the literature; see, e.g., \cite{Cai-Cot-Ger-2012}. The combination of continuation methods with shooting methods not only allows to preserve precision but also has been shown to be robust for problems of low-thrust orbital transfers \cite{doi:10.1080/10556788.2011.593625} and deep space exploration \cite{doi:10.2514/1.G005865}. 
Whereas, the solution obtained by the combination of shooting methods and continuation methods cannot be guaranteed to be at least locally optimal unless sufficient conditions for optimality are taken into account \cite{chen20161,chen2016optimality}.  

In contrast to indirect methods, direct methods transform  the original optimal control problem to a finite-dimensional nonlinear constrained optimization problem via a discretization of the control and the state functions on a time grid. Then, nonlinear programming (NLP) solvers can be used to solve the finite optimization problem by satisfying Karush–Kuhn–Tucker (KKT) conditions. Compared with  indirect methods, the direct methods are more flexible and robust to deal with large systems but less accurate. In addition, it should be noted that the KKT conditions are necessary but not sufficient  for optimality. Thus, the solution trajectories obtained from NLP solvers cannot be guaranteed to be at least locally optimal unless the finite optimization problem  is convex. In addition to the issue of optimality, direct methods are computationally demanding and should be used offline, as analyzed by Bonalli, H\'eriss\'e, and Tr\'elat in \cite{BONALLI2017482,8765338}.

According to the summarizations in the preceding two paragraphs, both indirect and direct methods suffer the issue of convergence. This prohibits the two classes of methods to be used onboard or in situ. In order to realize the real-time generation of optimal trajectories, researchers in industry and academia have made much effort. For example, because convex optimization problems can be addressed within a constant time, some lossless convexification approaches have been proposed in the literature to convert original optimal control problem into convex optimization problems \cite{ACIKMESE2011341}. Whereas, some nonlinear flight dynamics are hard to be losslessly convexified. Another way for real-time generation of optimal trajectories is to use Artificial Neural Networks (ANNs). In fact, as an ANN can generate its output within a constant time, some different architectures of ANNs have been proposed to generate optimal trajectories in aerospace engineering, such as orbital transfer of satellites \cite{izzo2021real}, powered descent guidance of rockets \cite{doi:10.2514/1.G004928,sanchez2018real}, trajectory optimization of hypersonic vehicles \cite{doi:10.2514/1.G004928}, just to name a few.


The common way of using ANN to generate the optimal solution of an optimal control problem is to train an ANN by a large number of optimal trajectories. Then, the trained ANN is embedded in onboard computers to generate optimal feedback control in real time. The dataset of optimal trajectories  are usually obtained by the aforementioned indirect or direct methods offline. The shooting methd, combined with homotopy method, was employed to generate sample data of optimal trajectories for hypersonic vehicles \cite{doi:10.2514/6.2020-0023}, orbital transfer \cite{cheng2018real}, and deep space exploration \cite{cheng2020real}. The direct method was used to generate dataset for powered landing \cite{wang2023real,sanchez2016learning} and midcourse guidance \cite{tang2021midcourse}. However, as stated above, both indirect and direct methods suffer the issue of convergence, making it time-consuming to generate the dataset for training ANNs. Recently, a backward propagation was proposed by Izzo and {\"O}zt{\"u}rk \cite{izzo2021real}, relying on integrating the canonical differential equations derived from the PMP in backward way. Although this method skillfully avoids the convergence issue  during generating the dataset for ANN, it cannot guarantee every solution in the dataset to be at least locally optimal because only necessary conditions are met.


Stacking none optimal solutions in the dataset prohibits using ANNs because the output of the trained network may not be correct. In addition, the training process may not converge.  In order to make sure that the trajectories for training ANNs are at least locally optimal, this paper proposes a parametrization approach so that not only necessary conditions but also sufficient conditions for optimality are embedded into a parameterized set of differential equations. By establishing the initial conditions for the parameterized differential equations, a simple numerical procedure is developed, allowing to generate the dataset of at least locally optimal trajectories via solving some simple initial value problems. Once a simple feedforward neural network is trained by the dataset, it can not only generate optimal trajectories within milliseconds but also ensure the generated trajectories to be at least locally optimal,  as finally demonstrated by two conventional MECPs in aerospace engineering.

The paper is organized as follows. In Section \ref{SE:Problem}, the MECP is formulated, and the basic principle for ANNs to generate optimal control command in real time is stated. Section \ref{SE:Parameterization} is devoted to establishing the parameterized family of extremals by embedding sufficient conditions for optimality, and Section \ref{SE:Procedure} presents how to use the parameterized extremals to generate the dataset for the mapping from flight state to the optimal feedback control. Two numerical examples are provided in Section \ref{SE:Numerical} to demonstrate in detail how to use the developed method to get real-time solutions for MECPs. 

\section{Preliminaries}\label{SE:Problem}

Throughout the paper, we denote the space of $n$-dimensional column vectors by $\mathbb{R}^n$, and the space of $n$-dimensional row vectors by $(\mathbb{R}^n)^*$.

\subsection{Optimal Control Problem}

Let $\mathcal{X}$ be a smooth manifold of dimension $n$. Consider the system
\begin{align}
(\Sigma):\ \ \ \ \dot{\x}(t) = \f(\x(t),\u(t))
\label{EQ:system}
\end{align}
where $t\in \mathbb{R}_+$ denotes time, the over dot denotes the differentiation with respect to time, $\x\in \mathcal{X}$ denotes the state,  $\u$ denotes the control vector, taking values in an open subset $\mathcal{U}$ of $\mathbb{R}^m$, and $\f:\mathcal{X}\times \mathcal{U}\rightarrow \mathbb{R}^n$ is a smooth vector field on $\mathcal{X}$. 
Let $s\leq n$ be a positive integer. Then, we define the constraining submanifold of the final state   by
\begin{align}
\mathcal{X}_f \triangleq \{\x \in \mathcal{X}| \phi(\x ) = 0\}
\end{align}
where $\phi:\mathcal{X}\rightarrow \mathbb{R}^{s}$ is a twice continuously differentiable function. Without loss of generality, we make a regular assumption on $\mathcal{X}_f$ as below.
\begin{assumption}
The matrix $\nabla \phi(\x)$ is of full rank for any $\x\in \mathcal{X}_f$, where the notation ``~$\nabla$'' is the gradient operator.
\label{AS:full_rank} 
\end{assumption}
Let $\x_c\in \mathcal{X}$ be the current state, and let $t_c > 0$ be the current time.  Then, we consider the following optimal control problem: 
\begin{problem}\label{Problem}
The MECP consists of finding a measurable control $\u(\cdot):[t_c,t_f]\rightarrow \mathcal{U}$ that steers $(\Sigma)$ from the current state $\x_c$ to a point $\x_f \in \mathcal{X}_f$ so that the control effort is minimized, i.e., 
\begin{align}
\begin{split}
J= \int_{t_c}^{t_f} \|\u(t)\| ^2 \mathrm{d}t \rightarrow \mathrm{min}
\end{split}
\end{align}
where $t_f \in (t_c, + \infty)$ is the specified final time and the notation $\| \cdot \|$ is the usual Euclidean norm. 
\end{problem}

\par The MECP is also called as minimum-energy control problem in the literature. As stated in Section \ref{SE:Introduction}, real-time solution of the MECP is required in many practical applications, but it is quite challenging to address the MECP within a constant time by onboard computers. In the following subsection, the principle of using ANN to generate real-time solutions of MECP will be presented.

\subsection{Principle for Real-Time Solutions via Neural Networks}

The value function for the MECP in Problem \ref{Problem} can be written as
\begin{align}
V(t,\boldsymbol{x}(t)) \triangleq \underset{\u\in \mathcal{U}}{\max}\int_{t}^{t_f} -\|\u(\tau)\|^2 \mathrm{d} \tau
\end{align}
Then, the Hamilton-Jacobi-Bellman (HJB) equation is given by
\begin{align}
-\frac{\partial V(t,\x)}{\partial t} =  \underset{\u\in \mathcal{U}}{\max}\Big \{ -\|\u\|^2 +\frac{\partial V(t,\x)}{\partial \x^T} \f(\x,\u)\Big\}
\label{EQ:HJB}
\end{align}
A solution of the HJB equation easily provides the optimal feedback control policy as the optimum of the inner maximization:
\begin{align}
\u^*(t,\x(t)) \triangleq \underset{\u \in \mathcal{U}}{\mathrm{argmax}}\Big\{ - \|\u\|^2 + \frac{\partial V(t,\x(t))}{\partial \x^T}\f(\x(t),\u)\Big\}
\label{EQ:u*_inner_minimization}
\end{align}

For the sake of notational simplicity, let $t_g\geq 0$ be the time to go, i.e., 
$$t_g = t_f - t_c$$
Then, for any current state $\x_c$, the time-to-go $t_g$ is said to be feasible if there exists an optimal trajectory $\x(\cdot):[0,t_f]\rightarrow \mathcal{X}$ of the MECP so that $\x_c = \x(t_f - t_g)$. Furthermore, we denote by $\mathcal{F}\subset \mathcal{X}\times[t_c,t_f]$ the set of $(t_g,\x_c)$ for which $t_g$ is feasible. Let $C(t_g,\x_c):\mathcal{F}\rightarrow \mathcal{U}$ be the optimal feedback control of the MECP at $(t_g,\x_c)$. Then, according to Eq.~(\ref{EQ:u*_inner_minimization}), we immediately have 
\begin{align}
C(t_g,\x_c) = \underset{\u \in \mathcal{U}}{\mathrm{argmax}}\Big\{-\|\u\|^2 + \frac{\partial V(t_f - t_g,\x_c)}{\partial \x^T}\f(\x_c,\u)\Big\}
\label{EQ:inner_minimization2}
\end{align}
Apparently, solving the MECP in real time is equivalent to addressing the maximization problem on the right side of Eq.~(\ref{EQ:inner_minimization2})  within a short constant time. However, it is well known that the gradient value $\partial V(t,\x)/\partial \x$ is hard to obtain as solving the nonlinear PDE in Eq.~(\ref{EQ:HJB}) is in general a very intractable task. 

\begin{figure}[!htp]
\centering
\includegraphics[width=0.4\textwidth]{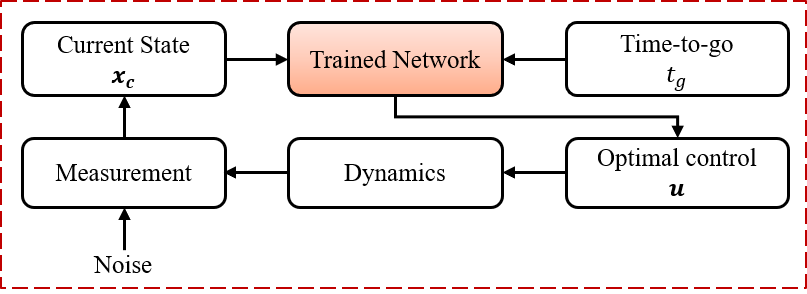}
\caption{Diagram for real-time solution via FNN.}
\label{Fig:closed_loop}
\end{figure}

If one is able to solve the inner maximization problem in Eq.~(\ref{EQ:inner_minimization2}) offline, one can obtain the dataset for the mapping from  $(t_g,\x_c)$ to the corresponding optimal feedback control $C(t_g,\x_c)$. Then, in view of the universal approximation theorem established in \cite{hornik1989multilayer,cybenko1989approximation,hornik1991approximation},  a simple Feedforward Neural Network (FNN) can be trained by the dataset to approximate the mapping $(t_g,\x_c)\mapsto C(t_g,\x_c)$ in Eq.~(\ref{EQ:inner_minimization2}). Notice that the output of an FNN is a composition of some linear mappings  of the input vector. Thus, for any feasible pair $(t_g,\x_c)\in \mathcal{F}$ as the input vector, the trained FNN can  generate the corresponding optimal feedback control  $C(t_g,\x_c)$ within a constant time, and it can be embedded in a closed-loop guidance system, as shown in Fig.~\ref{Fig:closed_loop}.  



As stated in Sect. \ref{SE:Introduction},  both indirect methods and direct methods have been used to generate the dataset. However, they cannot guarantee the solutions to be at least locally optimal. In the following two sections, an parametrization approach  will be developed so that not only necessary but also sufficient conditions are embedded into a parameterized set of differential equations. As a consequence, one is able to generate at least locally optimal trajectories of the MECP  via simply solving some initial value problems.


\section{Parametrization of Hamiltonian Extremals}\label{SE:Parameterization}






Note that the maximand on the right-hand side of Eq.~(\ref{EQ:HJB}) can be rewritten as the Hamiltonian
\begin{align}
H(\x,\p,\u) \triangleq \p \f(\x,\u) - \|\u\|^2
\label{EQ:Hamiltonian}
\end{align}
with 
\begin{align}
\p(t)=\partial V(t,\x)/\partial \x^T
\label{EQ:p}
\end{align}
 playing the role of the costate vector. 
By differentiating Eq.~(\ref{EQ:p}) with respect to time, we have 
\begin{align}
\dot{\p}(t) = \frac{\partial^2 V(t,\x)}{\partial \x^T \partial t} +\f^T(\x,\u) \frac{\partial^2 V(t,\x)}{\partial \x^T \partial \x}
\label{EQ:dot_pp}
\end{align}
Then, differentiating both sides of the HJB equation in Eq.~(\ref{EQ:HJB}) with respect to $\x$ leads to
\begin{align}
-\frac{\partial^2 V(t,\x)}{\partial t\partial \x} = \f^T(\x,\u) \frac{\partial^2 V(t,\x)}{\partial \x^T \partial \x}+\frac{\partial V(t,\x)}{\partial \x^T} \frac{\partial \f(\x,\u)}{\partial \x}
\label{EQ:dot_V}
\end{align}
Substituting Eq.~(\ref{EQ:dot_pp}) into Eq.~(\ref{EQ:dot_V}) and rewriting the system in Eq.~(\ref{EQ:system}), we have that each optimal trajectory $\x(t)$ and its corresponding optimal feedback control $\u^*(t,\x(t))$ for $t\in [0,t_f]$ satisfy the following differential equations:
\begin{align}
\begin{split}
\dot{{\x}}(t) =&\ \frac{\partial H}{\partial \p^T}({\x}(t),\p(t),\u^*(t,\x(t))),\\
\dot{{\p}}(t) =&\ -\frac{\partial H}{\partial \x^T}({\x}(t),{\p}(t),\u^*(t,\x(t))).
\end{split}
\label{EQ:cannonical}
\end{align}
Any pair $({\x}(\cdot),{\p}(\cdot))$  on $[0,t_f]$, solving the differential equations in Eq.~(\ref{EQ:cannonical}), is called a Hamiltonian extremal, and the corresponding control along each Hamiltonian extremal is said to be an extremal control.





Notice that the final boundary transversality condition, i.e., ${\boldsymbol{\p}}(t_f)\perp T_{{\x}(t_f)}\mathcal{X}_f$, implies
\begin{align}
{\boldsymbol{p}}(t_f) &=   {\boldsymbol{\nu}}_f \nabla \boldsymbol{\phi}_f({\x}(t_f)),\label{EQ:Transversality_2}
\end{align}
where ${\boldsymbol{\nu}}_f\in(\mathbb{R}^{s})^*$ is a constant vector whose elements are Lagrangian multipliers. 
Therefore, the submanifold 
\begin{align}
\mathcal{L}_f =\{(\x,\p)\in T^*\mathcal{X}\ \arrowvert\ \x\in\mathcal{X}_f,\ \p\perp T_{\x}\mathcal{X}_f\}
\label{EQ:L_f}
\end{align}
is a Lagrangian manifold in the cotangent bundle $T^*\mathcal{X}$.
Given any $(\x_f,\p_f) \in \mathcal{L}_f$, let
$$\gamma(t,\x_f,\p_f):[0,t_f]\times \mathcal{L}_f \rightarrow T^*\mathcal{X}$$
be the solution trajectory
\begin{align}
\begin{cases}
\dot{\x}(t) =\  - \frac{\partial H}{\partial \p}(\x(t),\p(t),\u(t,\x(t))),\\
\dot{\p}(t) =\  \frac{\partial H}{\partial \x}(\x(t),\p(t),\u(t,\x(t)))
\end{cases}
\label{EQ:Para_system}
\end{align}
 with the initial condition 
\begin{align}
\gamma(0,\x_f,\p_f)= (\x_f,\p_f)
\end{align}
Notice that the right side of the differential equations in Eq.~(\ref{EQ:Para_system}) is opposite to that in Eq.~(\ref{EQ:cannonical}). Thus, we have that $\gamma(t_f-t,\x_f,\p_f)$ for $t\in [0,t_f]$ is a Hamiltonian extremal for every $(\x_f,\p_f)\in  \mathcal{L}_f$. Given any $(\x_f,\p_f)\in \mathcal{L}_f$, let us denote by $T_c(\x_f,\p_f)$ the conjugate time or focal time. That means the extremal trajectory $\gamma(t,\x_f,\p_f)$ for $t\in (t_f-T_c(\x_f,\p_f),t_f]$ is at least locally optimal. Then, by the definition of conjugate time, we have the following conclusions:
\begin{remark}[Chen \cite{chen2017neighboring}]\label{RE:optimal}
Given any pair $(\x_f,\p_f)$ in $\mathcal{L}_f$, if $T \leq T_c(\x_f,\p_f)$, the Hamiltonian extremal  $\gamma(t_f-t,\x_f,\p_f)$  for $t\in [0,T]$ is locally optimal; however, it loses its local optimum if $T > T_c(\x_f,\p_f)$.
\end{remark}
Set 
\begin{align}
\begin{split}
T({\x}_f,{\p}_f)
\triangleq  \min \Big\{ t_f, T_c(\x_f,\p_f) \Big\}
\end{split}
\label{EQ:T}
\end{align}
and set
\begin{align}
\begin{split}
\mathcal{L} \triangleq \big\{\gamma(t,\x_f,\p_f) \in  T^*\mathcal{X}\ \arrowvert t\in (0,T(\x_f,\p_f)],(\x_f,\p_f)\in \mathcal{L}_f\big\}
\end{split}
\label{EQ:graph}
\end{align}
Up to now, all the Hamiltonian extremals have been embedded into the parameterized family $\mathcal{L}$. According to the conclusions in Remark \ref{RE:optimal}, not only necessary conditions but also sufficient conditions are met along each extremal trajectory in the family $\mathcal{L}$. In the next section, we shall show how to generate the dataset for the mapping $(t_g,\x_c)\mapsto C(t_g,\x_c)$ via using the parameterized extremals in $\mathcal{L}$.

\section{Procedure for Generating the Dataset}\label{SE:Procedure}

Given any extremal trajectory in $\mathcal{L}$, we can immediately obtain some data for the mapping $(t_g,\x_c)\mapsto C(t_g,\x_c)$ by discretizing the extremal trajectory. Therefore, in order to generate the dataset for the mapping $(t_g,\x_c)\mapsto C(t_g,\x_c)$, it is enough to generate some extremal trajectories in $\mathcal{L}$. According to the developments in Section \ref{SE:Parameterization}, if we choose a pair $(\x_f,\p_f)\in \mathcal{L}_f$ as the initial condition, we can obtain an extremal trajectory by integrating the set of differential equations in Eq.~(\ref{EQ:Para_system}) from $t=0$ to $t=T(\x_f,\p_f)$. Whereas, the value of $T(\x_f,\p_f)$ is not trivial to obtain as it relies on computing the conjugate time $T_c(\x_f,\p_f)$. In this section, the numerical method for calculating the conjugate time will first be presented, which will allow us to develop a propagation method for generating the dataset of the mapping $(t_g,\x_c)\mapsto C(t_g,\x_c)$.

\subsection{Calculation of Conjugate Time}

Without loss of generality, let us arbitrarily choose
 a pair $(\bar{\x}_f,\bar{\p}_f)\in \mathcal{L}_f$.
 Due to Assumption \ref{AS:full_rank}, there exists a local coordinate chart $\boldsymbol{F}$ of $\mathcal{L}_f$ around $(\bar{\x}_f,\bar{\p}_f)$. Then, for any $(\x_f,\p_f)$ in a sufficiently small neighborhood of $( \bar{\x}_f, \bar{\p}_f)$ there exists one and only one $\q \in \mathbb{R}^n$ satisfying $\q = \boldsymbol{F}(\x_f,\p_f)$. Without lose of generality, let us choose a coordinate chart $\boldsymbol{F}$ so that $\boldsymbol{F}(\bar{\x}_f,\bar{\p}_f) = \boldsymbol{0}$. Set
\begin{align}
[X(t,\q),P(t,\q)] = \gamma(t,\boldsymbol{F}^{-1}(\q))
\end{align}
Apparently, for any $\q\in \mathbb{R}^n\setminus\{\boldsymbol{0}\}$ the pair $[X(t,\q),P(t,\q)]$  is a Hamiltonian extremal in $\mathcal{L}$. According to \cite{chen2017neighboring}, the conjugate point occurs if the matrix $\partial X(t,\boldsymbol{0})/\partial \q $ losses its full rank. Thus, we can calculate the conjugate time $T_c(\bar{\x}_f,\bar{\p}_f)$ by finding when the determinant of  $\partial X(t,\boldsymbol{0})/\partial \q $ is zero.




By the definition of the parameterized Hamiltonian extremal $[X(t,\q),P(t,\q)]$ and according to Eq.~(\ref{EQ:Para_system}), we have
\begin{align}
\begin{split}
\dot{X}(t,\q) & = - \frac{\partial H}{\partial \p^T}[X(t,\q),P(t,\q),\u(t,X(t,\q))]\\
\dot{P}(t,\q) &= \frac{\partial H}{\partial \x^T}[X(t,\q),P(t,\q),\u(t,X(t,\q))]
\end{split}
\label{EQ:Para_system1}
\end{align}
It follows from the classical results about solutions to ODEs that the Hamiltonian extremal  $({X}(\cdot,{\q}),{P}(\cdot,{\q}))$ and its time derivative are continuously differentiable with respect to $\q$. Thus, taking the derivative of Eq.~(\ref{EQ:Para_system1}) with respect to $\q$, we obtain the homogeneous linear matrix differential equations
\begin{align}
\begin{split}
\frac{\mathrm{d}}{\mathrm{d}t} \frac{\partial {X}}{\partial \q} (t,\boldsymbol{0}) = &\ - H_{\p\x}(\bar{\x}(t),\bar{\p}(t)) \frac{\partial  {X}}{\partial \q} (t,\boldsymbol{0})\\
- &\  H_{\p\p}(\bar{\x}(t),\bar{\p}(t))\frac{\partial  {P}}{\partial \q} (t,\boldsymbol{0})\\
\frac{\mathrm{d}}{\mathrm{d}t} \frac{\partial  {P}}{\partial \q} (t,\boldsymbol{0}) =  &\ H_{\x\x}(\bar{\x}(t),\bar{\p}(t))  \frac{\partial  {X}}{\partial \q} (t,\boldsymbol{0})\\
+ &\ H_{\x\p}(\bar{\x}(t),\bar{\p}(t))\frac{\partial {P}}{\partial \q} (t,\boldsymbol{0})
\end{split}
\label{EQ:matrix2}
\end{align}
where $\bar{\x}(t) = X(t,\boldsymbol{0})$ and $\bar{\p}(t) = P(t,\boldsymbol{0})$.
Once the initial conditions  $\partial X(0,\boldsymbol{0})/\partial \q$ and $\partial P(0,\boldsymbol{0})/\partial \q$ are given, one can obtain the two matrices $\partial X(t,\boldsymbol{0})/\partial \q$ and $\partial P(t,\boldsymbol{0})/\partial \q$ for any time $t\in [0,t_f]$ by integrating the differential equations in Eq.~(\ref{EQ:matrix2}).

\subsubsection{Initial Conditions for $s= n$}

 In order to get the initial conditions $\partial X(0,\boldsymbol{0})/\partial \q$ and $\partial P(0,\boldsymbol{0})/\partial \q$, it is enough to compute a set of basis vectors of the tangent space of the manifold $\mathcal{L}_f$ at $(\bar{\x}(t_f),\bar{\p}(t_f))$. In view of Eq.~(\ref{AS:full_rank}), the submanifold $\mathcal{X}_f$ reduces to a singleton if $s = n$.  Hence, in the case of $s = n$, one can simply set $\q = \p_f - \bar{\p}_f$, which implies 
\begin{align}
\frac{\partial X(0,\boldsymbol{0})}{\partial \q} = O_n\ \text{and}\ \frac{\partial P(0,\boldsymbol{0})}{\partial \q} = I_n
\label{EQ:matrix_initial}
\end{align}
where $O_n$ and $I_n$ denote the zero and identity matrices of $\mathbb{R}^{n\times n}$, respectively. 

\subsubsection{Initial Conditions for $s< n$}

If $s < n$, there exists an invertible function $\hat{\boldsymbol{F}}:\mathcal{X}_f\rightarrow (\mathbb{R}^{n-s})^*$ so that both the function and its inverse $\hat{\boldsymbol{F}}^{-1}$ are smooth. Note that the function $\hat{\boldsymbol{F}}$ is a coordinate chart on $\mathcal{X}_f$. Then, for every $\x_f$ in a small neighborhood of $\bar{\x}_f$ there exists one and only one $\hat{\q}$ so that $\hat{\q}  = \hat{\boldsymbol{F}}(\x_f)$. According to the transversality condition in Eq.~(\ref{EQ:Transversality_2}), for every $(\x_f,\p_f)\in \mathcal{L}_f$,  there exists a $\boldsymbol{\nu}\in (\mathbb{R}^{s})^*$ so that 
\begin{align}
\p_f = \boldsymbol{\nu}\nabla \phi(\x_f)
\end{align}
Then, it is enough to set $\q = [\hat{\q} - \hat{\boldsymbol{F}}(\bar{\x}_f),\boldsymbol{\nu}-\bar{\boldsymbol{\nu}}]$, where
\begin{align}
\bar{\boldsymbol{\nu}} = \bar{\p}_f \nabla \phi^T(\bar{\x}_f) [\nabla \phi(\bar{\x}_f)\nabla^T (\bar{\x}_f)]^{-1}
\end{align}
denotes the vector of the Lagrangian multipliers in Eq.~(\ref{EQ:Transversality_2}) for the final point $(\bar{\x}_f,\bar{\p}_f)$. Then, a direct calculation leads to
\begin{align}
\frac{\partial X}{\partial \q}(0,\boldsymbol{0}) = &\ \left[\frac{\partial X}{\partial \hat{\q}}(0,\boldsymbol{0}),\frac{\partial X}{\partial \boldsymbol{\nu}}(0,\boldsymbol{0})\right]\label{EQ:final_condition1}\\
\frac{\partial P}{\partial \q}(0,\boldsymbol{0}) =&\  \left[\frac{\partial P^T}{\partial \hat{\q}}(0,\boldsymbol{0}),\frac{\partial P^T}{\partial \boldsymbol{\nu}}(0,\boldsymbol{0})\right]\nonumber\\
= &\ \left[\bar{\boldsymbol{\nu}}\nabla^2\phi(\bar{\x}_f)\frac{\partial X}{\partial \hat{\q}}(0,\boldsymbol{0}),\nabla\phi^T(\bar{\x}_f)\right]
\label{EQ:final_condition2}
\end{align}
Since $X(0,\q)$ is not a function of $\boldsymbol{\nu}$, it follows that 
\begin{align}
\frac{\partial X}{\partial \boldsymbol{\nu}}(0,\boldsymbol{0}) = O_{n\times  s}
\label{EQ:matrix_initial1}
\end{align}
where $O_{n\times s}$ denotes the zero matrix in $\mathbb{R}^{n\times s}$. Up to present, all the quantities for computing the conditions in Eq.~(\ref{EQ:final_condition1}) and Eq.~(\ref{EQ:final_condition2}) are available except the matrix $\partial X(0,\boldsymbol{0})/\partial \hat{\q}$.  Let us take the differentiation of $\phi(X(0,\q)) = 0$ with respect to $\hat{\q}$. Then, we have
\begin{align}
\nabla \phi(X(0,\boldsymbol{0}))\frac{\partial X(0,\boldsymbol{0})}{\partial \hat{\q}} = 0
\label{EQ:matrix_initial2}
\end{align}
This equation implies that all the column vectors of the matrix $\partial X(0,\boldsymbol{0})/\partial \hat{\q}$ constitutes a basis of the tangent space $T_{\bar{\x}_f}\mathcal{X}_f$. Once the matrix $\nabla \phi(X(0,\boldsymbol{0}))$ is given, one can compute the full-rank matrix $\partial X(0,\boldsymbol{0})/\partial \hat{\q}$ by a Gram-Schmidt orthogonalization.

\subsection{Numerical Procedure for Generating the Dataset}


According to the developments in the previous subsection, given any $(\x_f,\p_f)\in \mathcal{L}_f$, we are able to obtain an optimal trajectory via solving an initial value problem. As a consequence, by sampling some pairs in $\mathcal{L}_f$ as initial conditions, we are able to use the initial value problem to generate sampled data for Hamiltonian extremals in $\mathcal{L}$. Then, the dataset for the mapping from the flight state to the corresponding optimal feedback control can be immediately obtained, as shown in {\it Procedure \ref{algo1}}.
\usetikzlibrary{calc}
\begin{center}
	\begin{tcolorbox}[
		colframe=blue!25,
		colback=gray!20,
		coltitle=blue!20!black,  
		fonttitle=\bfseries,
		width = 1\columnwidth,
		boxrule=0pt,
		top=1pt,
		bottom=-0.5pt,
		segmentation code={\draw[black,solid,line width=1.2pt]($(segmentation.west)+(0.2,0)$)--($(segmentation.east)+(-0.2,0)$);}]
		\begin{example}
		\label{algo1}
		\end{example}
\begin{center}
\textbf{Procedure 1}: {\it Generation of the Dataset}
\end{center}
\bigskip

\begin{itemize}
\item[1.] Uniformly choose $N$ points from $\mathcal{L}_f$, and let us denote the chosen points by $(\x_f^i,\p_f^i)$, $i=1,2,\ldots,N$. Let $\Delta t$ be a positive number.
\item[2.] Set $i=1$ and $\mathcal{D}= \varnothing$. 
\item[3.] If $i\leq N$, go to step 4; otherwise, go to step 7.
\item[4.] Propagate the system in Eq.~(\ref{EQ:Para_system1}) from the initial condition $(\x_f^i,\p_f^i)$ over the interval $[0,T(\x_f^i,\p_f^i)]$, to generate the Hamiltonian extremal $(\x,\p)$ and corresponding extremal control $\u$. Set $t=0$ and go to step 5.
\item[5.] If $t+\Delta t \leq T(\x_f^i,\p_f^i)$, set $t = t + \Delta t $ and go to step 6; otherwise, set $i=i+1$ and go to step 3.
\item[6.] Set $\mathcal{D} = \mathcal{D}\cup \{[t,\x(t),\u(t)]\}$  and go to step 5.
\item[7.] End.
\end{itemize}
	\end{tcolorbox}
\end{center}
By {\it Procedure} \ref{algo1}, the dataset is finally included in set $\mathcal{D}$. A simple FNN trained by the dataset in $\mathcal{D}$ is able to approximate the mapping  $(t_g,\x_c)\rightarrow C(t_g,\x_c)$, as shown by the numerical examples in the next section.



\section{Numerical Applications}\label{SE:Numerical}

In this section, the theoretical developments will be demonstrated by two conventional MECPs in aerospace engineering, including optimal gliding of flight vehicles \cite{doi:10.2514/6.1998-4114} and  optimal proximity of satellites \cite{doi:10.2514/1.14580,Topputo:Approximate:2013}. 

\subsection{Application to Optimal Gliding of Flight Vehicles}

\begin{figure}[htbp]
\centering
\includegraphics[width=0.4\textwidth]{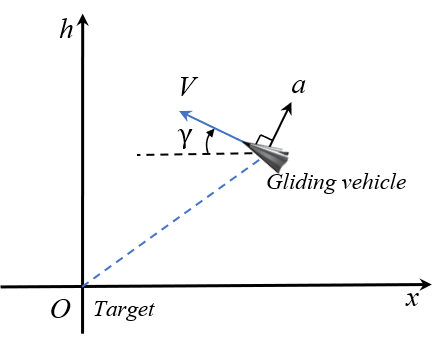}
\caption{Coordinate system for gliding of flight vehicle}
\label{Fig:Coordinate}
\end{figure}

We consider a scenario that a flight vehicle glides in a vertical plane, as shown in Fig.~\ref{Fig:Coordinate}.
The differential equations for the motion of the  vehicle are given as
\begin{align}
\begin{split}
\dot{V} = &\ - \frac{ D}{m} - g \sin \gamma,\\
\dot{\gamma} = &\  \frac{ a -  g \cos \gamma}{V},\\
\dot{x} = & \ V \cos \gamma ,\\
\dot{h} = &\ V \sin \gamma , \\
\end{split}
\label{EQ:System}
\end{align}
where $x$ is the downrange to the origin, $h>0$ is the altitude, $\gamma\in [-\pi/2,\pi/2]$ is the flight path angle,  $V\geq 0$ is the speed, and $a$ is normal acceleration representing the control parameter. $m = 100$ kg and $g=9.8$ m/s$^2$ are the mass  and the Earth gravitational acceleration, respectively. The expression of drag $D$ is given as \cite{doi:10.2514/6.1998-4114}
\begin{align}
D = k_1 V^2 + k_2 \frac{a^2 }{V^2}
\end{align}
with 
\begin{align}
k_1 = &\ \frac{1}{2}\rho S C_{D_0},\\
k_2 = &\ 2\frac{k_m m^2}{\rho S}
\end{align}
where $S=0.0324$ m$^2$ is the reference area, $C_{D_0}=0.2$ and $k_m=0.1$ are the zero-lift drag and the induced drag factor, respectively.  
We consider to control the vehicle from an initial state to the origin with free final speed and free final flight path angle while minimizing the control effort, i.e., 
\begin{align}
	\mathrm{min} \rightarrow \int_0^{t_f} a^2 dt
\end{align}
where $t_f$ is the expected flight time for the vehicle to the origin. 

The problem described above is exactly an MECP. Thus, the parameterized method summarized in {\it Procedure \ref{algo1}} can be used to generate the the dataset $\mathcal{D}$ for the mapping of optimal feedback control: $(t_g,\x_c)\mapsto C(t_g,\x_c)$. To this end, Set the integer $N$ and $\Delta t$ in {\it Procedure \ref{algo1}} as $6.2\times10^{5}$ and $0.5$ sec, respectively. This means that we uniformly choose  $6.2\times10^{5}$ points in $\mathcal{L}_f$. Then, the dataset $\mathcal{D}$ is directly generated by {\it Procedure \ref{algo1}}. An FNN with three hidden layers (each of which contains $20$ neurons) is trained by dataset $\mathcal{D}$ to approximate the optimal feedback control $C(t_g,\x_c)$. The training is terminated when the mean-squared error between the predicted values and the values in $\mathcal{D}$ is less than $10^{-6}$. Given any feasible pair $(t_g,\x_c)\in \mathcal{F}$ as input, the trained FNN takes around $0.15$ ms to generate an output on an embedded system with MYC-Y6ULY2 CPU at 528 MHz. Notice that such a computational period is more than enough for the onboard GNC system of usual flight vehicles.



\begin{figure}[!htp]
	\centering
	\includegraphics[width=0.8\textwidth]{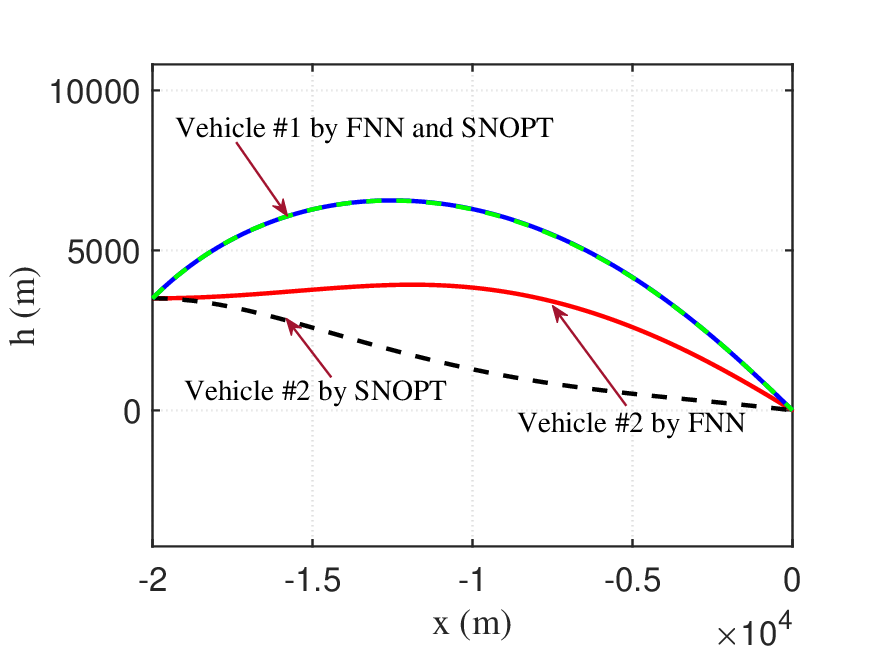}
	\caption{Trajectories for gliding of flight vehicles.}
	\label{Fig:salvo_trajectories}
\end{figure}
\begin{figure}[!htp]
	\centering
	\includegraphics[width=0.8\textwidth]{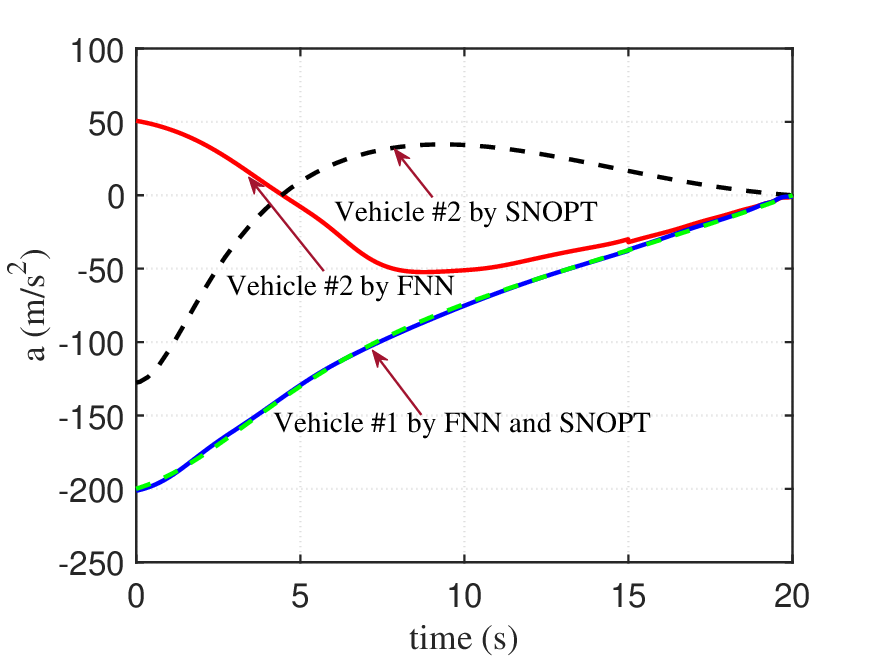}
	\caption{Control profiles for gliding of flight vehicles.}
	\label{Fig:salvo_controls}
\end{figure}
\begin{figure}[!htp]
	\centering
	\includegraphics[width=0.8\textwidth]{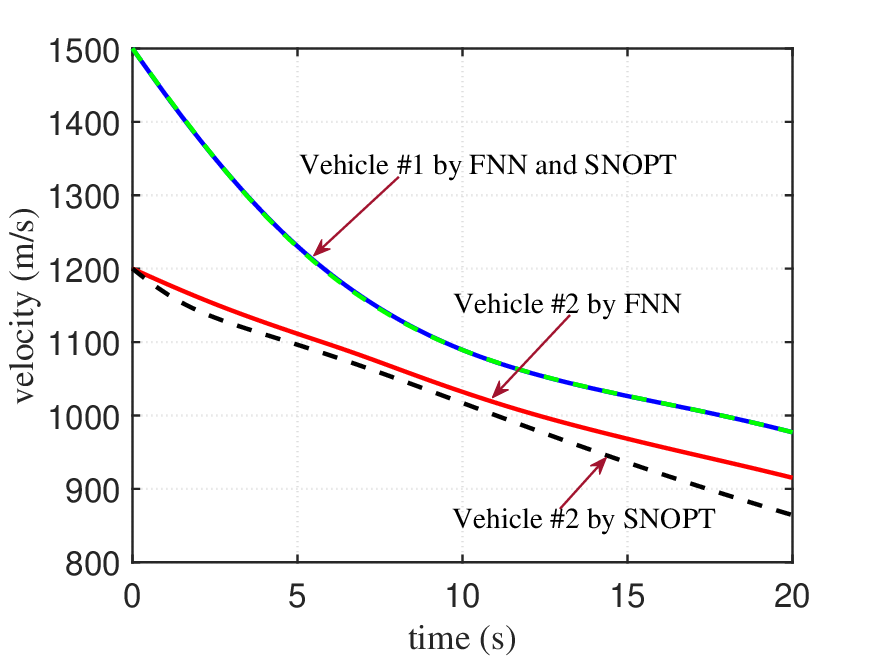}
	\caption{Time histories of speed for gliding of flight vehicles.}
	\label{Fig:salvo_speed}
\end{figure}
\begin{figure}[!htp]
	\centering
	\includegraphics[width=0.8\textwidth]{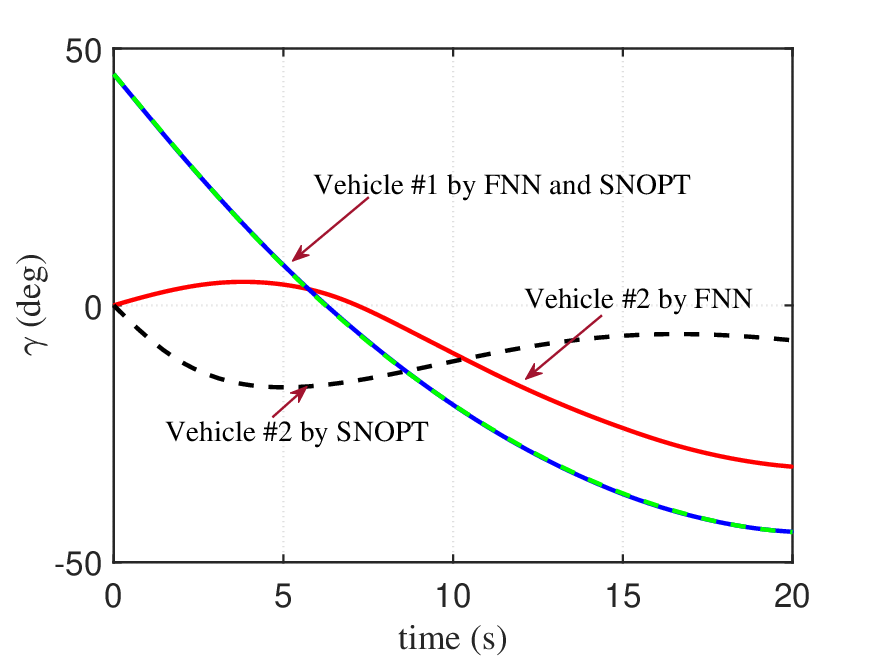}
	\caption{Time histories of flight path angle for gliding of flight vehicles.}
	\label{Fig:salvo_angle}
\end{figure}

We consider two different initial conditions, as presented in Table \ref{Tab:initial_condition}, to illustrate the trained FNN by comparing with optimization methods (here we employed the optimization toolbox of SNOPT).  
\begin{table}[!htp]
\centering
\caption{The initial states and initial time-to-go for each vehicle.}
\begin{threeparttable}[b]
	\label{Tab:initial_condition}
	\begin{tabular}{cccccc}
		\toprule
		States & $x_0$  & $h_0$ & $V_0$  & $\gamma_0$ & $t_{g0}$\footnotemark[1] \\
		\hline
		Vehicle $\# 1$ & $-20$ km & $3.5$ km  & $1500$ m/s &  $45$ deg  & 20 sec \\
		Vehicle $\# 2$ & $-20$ km & $3.5$ km  & $1200$ m/s &  $0$ deg  & 20 sec  \\
		\bottomrule
	\end{tabular}
	\begin{tablenotes}
		\item[1] $t_{g0}$ denotes the initial time-to-go.
	\end{tablenotes}
\end{threeparttable}
\end{table}
Note that the trained FNN can only generate the optimal feedback control command instead of the whole optimal trajectory. Thus, we embed the trained FNN into the closed-loop diagram in Fig.~\ref{Fig:closed_loop} in order to use the trained FNN to generate the whole optimal trajectory. The trajectories  generated by SNOPT (dashed curves) and trajectories generated by FNN (solid curves) are presented in Fig.~\ref{Fig:salvo_trajectories}, and the corresponding control profiles are reported in Fig.~\ref{Fig:salvo_controls}. The time histories of speed and flight path angle are presented in Fig.~\ref{Fig:salvo_speed} and Fig.~\ref{Fig:salvo_angle}, respectively.

\begin{figure}[!htp]
	\centering
	\includegraphics[width=0.8\textwidth]{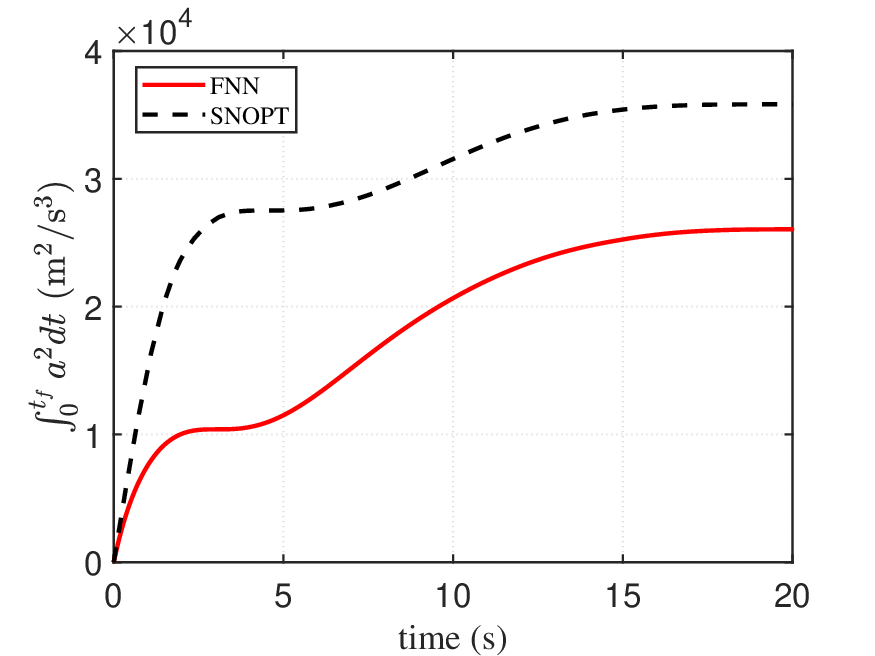}
	\caption{Control efforts for gliding of flight vehicle $\#2$.}
	\label{Fig:conj_u}
\end{figure}

It is apparent to see from Fig.~\ref{Fig:salvo_trajectories} that,  for vehicle $\#1$, the trajectories by FNN and SNOPT are the same; however, for vehicle $\#2$ the trajectory by FNN is quite different from that by SNOPT. To further examine which trajectory for vehicle $\#2$ in Fig.~\ref{Fig:salvo_trajectories} is optimal, the values of objective functions (control efforts) are presented in Fig.~\ref{Fig:conj_u}. It is seen from Fig.~\ref{Fig:conj_u} that the total control effort by SNOPT is $3.5830\times10^4$. However, the corresponding control effort by FNN is just $2.6048\times10^4$. We can also see from Fig.~\ref{Fig:salvo_controls} that the absolute value of normal acceleration by FNN is smaller than that by SNOPT. Therefore, from the perspective of expending less control efforts, the trajectories by FNN are better than SNOPT. 

In fact, the SNOPT belongs to the usual NLP methods, which generally satisfying KKT conditions (necessary conditions). Different from the NLP methods, the parameterized approach proposed in the current paper takes into account not only necessary conditions but also sufficient conditions, as guaranteed by the parameterized family $\mathcal{L}$ of extremal trajectories in Eq.~(\ref{EQ:graph}). It is also worth mentioning that even if the trajectory by SNOPT sometimes is the same as that by the trained FNN, as shown by the trajectories for vehicle $\#1$ in Fig.~\ref{Fig:salvo_trajectories}, the computation time for the two methods are quite different. The trained FNN can generate each optimal trajectory within a constant time. However, the SNOPT sometimes does not converge. 

\begin{figure}[!htp]
	\centering
	\includegraphics[width=0.8\textwidth]{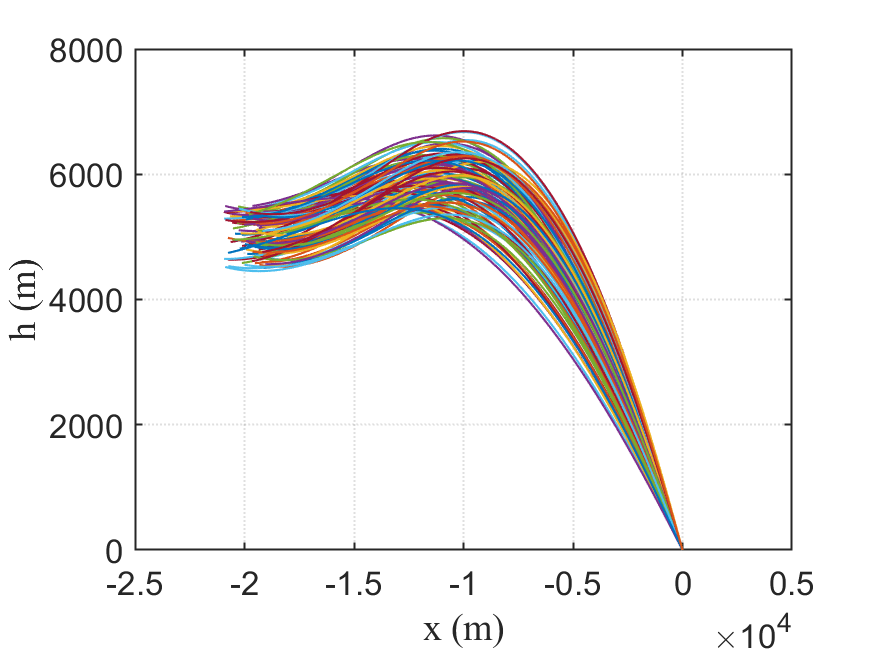}
	\caption{Trajectories of Monte Carlo test $\#1$.}
	\label{Fig:monte_trajectories1}
\end{figure}

\begin{figure}[htbp]
	\centering	
	\subfigure[Errors of terminal downrange.]{
		\includegraphics[width=0.8\textwidth]{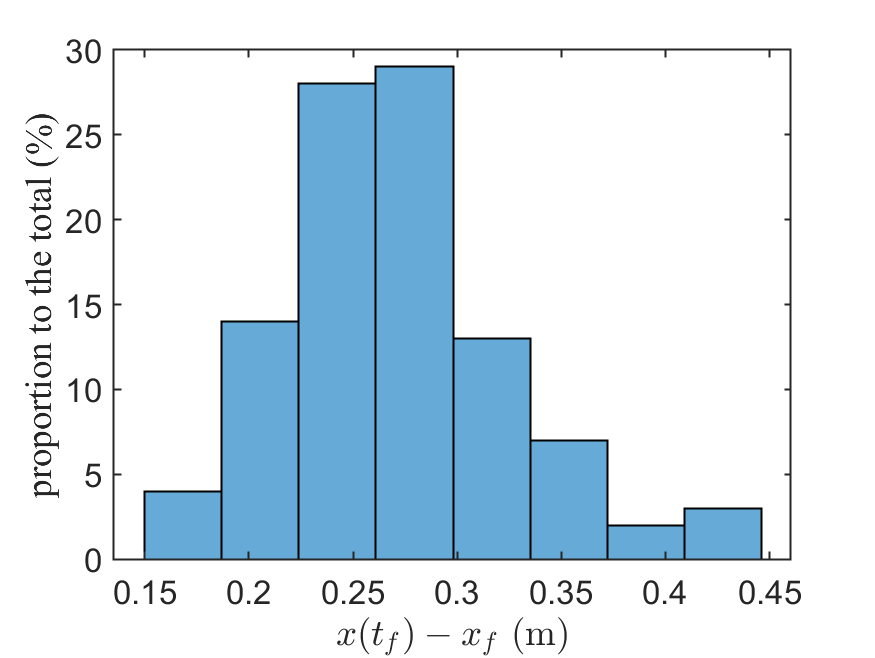}
		\label{Fig:monte_xerror}
	}
	\quad
	\subfigure[Errors between initial time-to-go and final time.]{
		\includegraphics[width=0.8\textwidth]{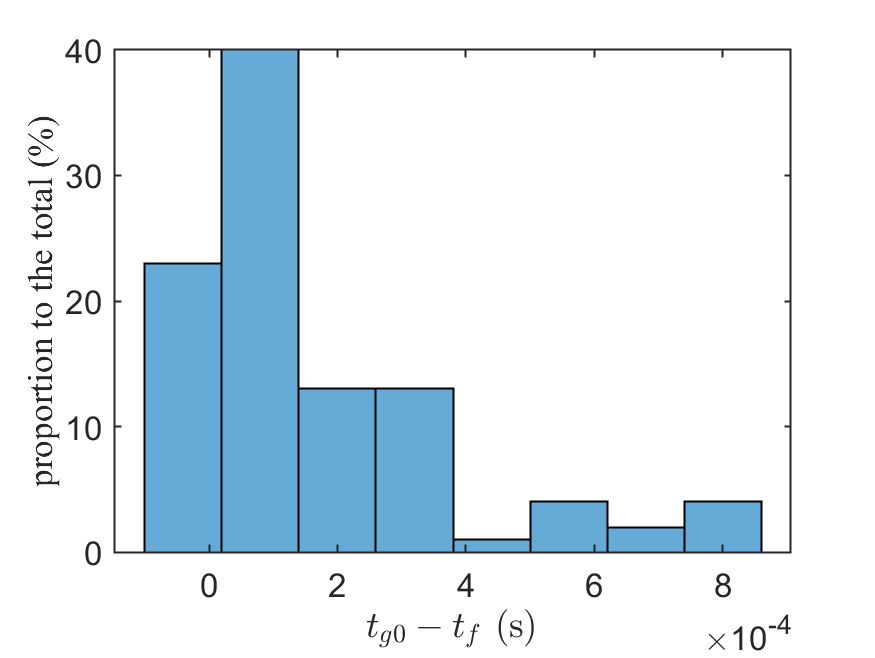}
	}
	\caption{Histograms for errors of Monte Carlo test $\#1$.}
	\label{Fig:monte_error}
\end{figure}

In order to demonstrate the robustness of the trained FNN, two Monte Carlo tests are implemented by disturbing the initial conditions and aerodynamic parameters. For test $\#1$, $100$ initial conditions are randomly selected. To be more specific, the initial speed, initial flight path angle, initial downrange, and initial altitude are randomly selected within the intervals $[1250, 1350]$ m/s, $[-0.1,0.1]$ rad, $[19, 21]$ km, and $[4,6]$ km, respectively. The trajectories generated by test $\#1$ are shown in Fig.~\ref{Fig:monte_trajectories1}, and the distribution of terminal errors are presented in Fig.~\ref{Fig:monte_error}. We can see from Fig.~\ref{Fig:monte_xerror} that the maximum error of downrange is less than $0.45$ m, and that the maximum error of expected flight time is less than $1$ ms. 

\begin{figure}[!htp]
	\centering
	\includegraphics[width=0.8\textwidth]{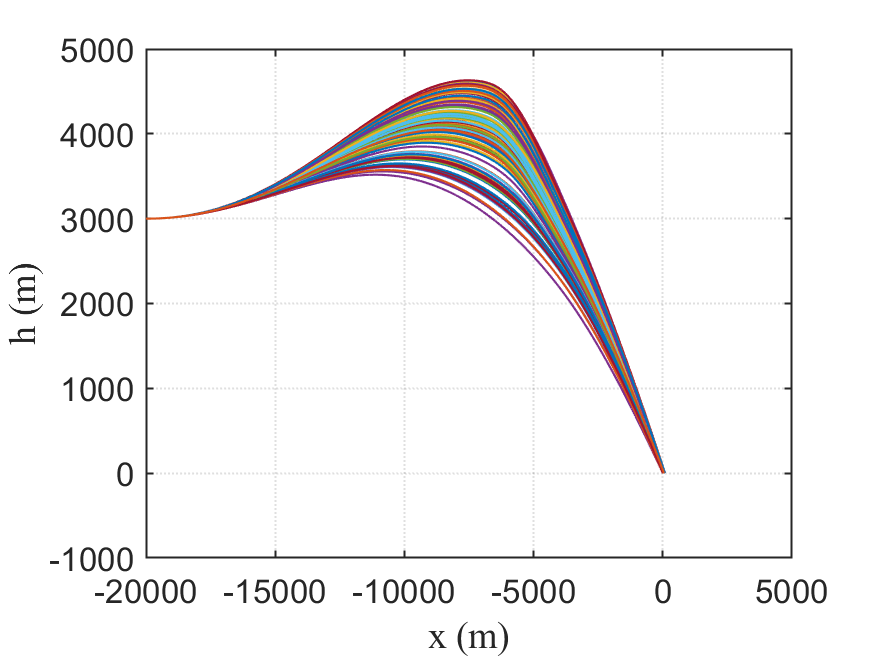}
	\caption{Trajectories of Monte Carlo test $\#2$.}
	\label{Fig:monte_trajectories2}
\end{figure}

\begin{figure}[htbp]
	\centering
	\subfigure[Errors of terminal downrange.]{
		\includegraphics[width=0.8\textwidth]{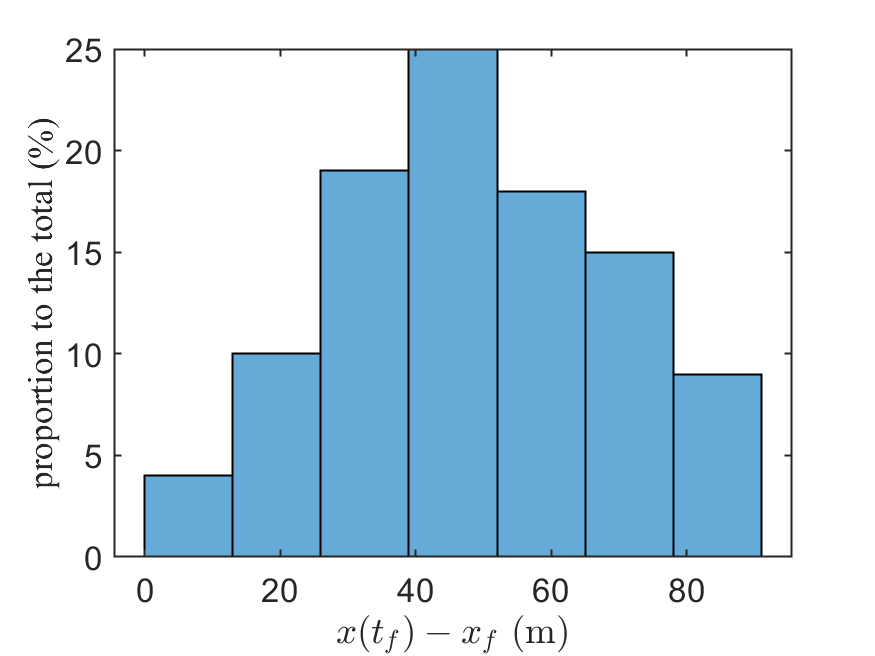}
	}
	\quad
	\subfigure[Errors between initial time-to-go and final time.]{
		\includegraphics[width=0.8\textwidth]{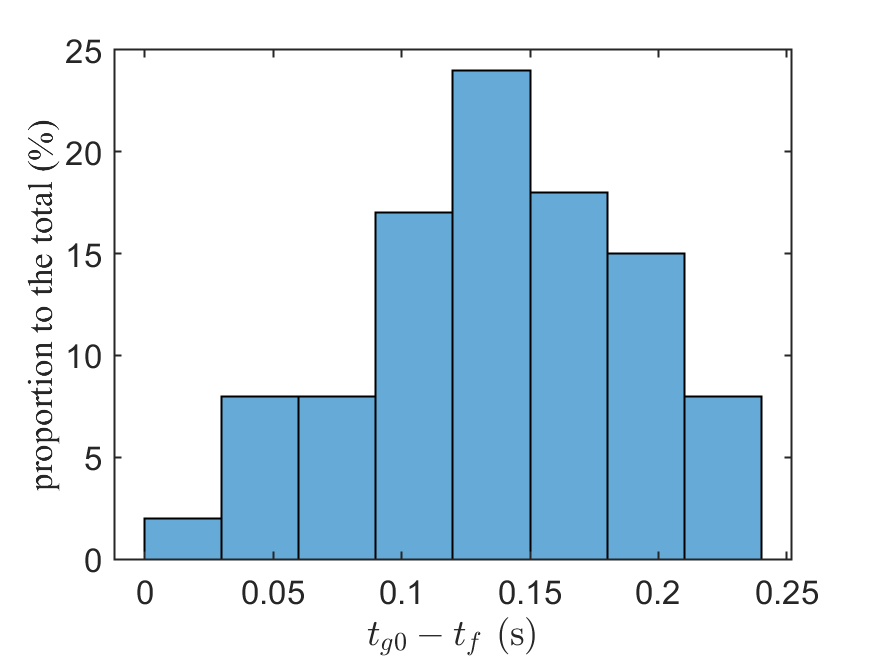}
	}
	\caption{Histograms for errors of Monte Carlo test $\#2$.}
	\label{Fig:monte_error0}
\end{figure}

Test $\#2$ is designed to show the robustness of the trained FNN by taking into consideration the disturbances of aerodynamics parameters. We consider that the initial state is fixed, but the values for the zero-lift drag $C_{D0}$ and the induced drag factor $k_m$ are randomly selected in $[0.1,0.2]$ and $[0.05,0.1]$, respectively. The trajectories of test $\#2$ with $100$ simulations are depicted in Fig.~\ref{Fig:monte_trajectories2}. And, the distribution of terminal errors is demonstrated by the histograms in Fig.~\ref{Fig:monte_error0}, from which we can see that the maximum error for the terminal downrange is less than $100$ m and the maximum error of expected flight time is less than $0.25$ s. 

\subsection{Application to Optimal Proximity of Spacecraft}

We consider a scenario of removing space debris on a circular orbit around the Earth, for which a crucial problem is to control one spacecraft, in a central gravity field, to a target (debris). For simplicity, the distance is normalized by the radius of the target's circular orbit. Thus, in the normalized setting, the period of the target's orbit is  $2\pi$, and the gravitational parameter is equal to 1. 

\begin{figure}[!htp]
	\centering
	\includegraphics[width=0.3\textwidth]{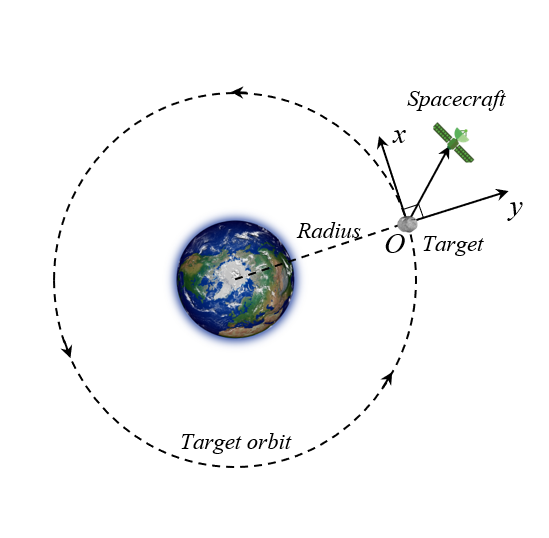}
	\caption{Coordinate system and geometry of spacecraft}
	\label{Fig:Coordinate2}
\end{figure}

Let us consider the local vertical local horizontal frame for which the $x$-axis is align with the target's velocity, and the $y$-axis points outward along the radius of the target, as presented in Fig.~\ref{Fig:Coordinate2}. Then, the equations of planar motion can be expressed as
\begin{align}
	\left[
	\begin{array}{c}
		\dot{x}\\
		\dot{y}\\
		\dot{v}_x\\
		\dot{v}_y
	\end{array}
	\right]
	=\left[
	\begin{array}{c}
		v_x\\
		v_y\\
		2v_y - (1-x)(1/r^3 - 1) + u_x\\
		-2 v_x - y(1/r^3 - 1) + u_y
	\end{array}\right]
\end{align}
where  $x$ and $y$ represent the radial and tangential displacements, respectively;  $v_x$ and $v_y$ denote the radial and tangential speed deviations, respectively; $u_x$ and $u_y$ are the radial and tangential accelerations, respectively; $r = \sqrt{(x+1)^2 + y^2}$ denotes the normalized distance to the center of the Earth. The performance index for the corresponding MECP is  
\begin{align}
	J=\frac{1}{2}\int_0^{t_f} u_x^2+u_y^2 \mathrm{d}t
\end{align}
In order to remove the debris, the final boundary condition has to be given by 
$(x(t_f),y(t_f)) = (0,0)$, and we set the normalized final time as 1, i.e., $t_f = 1$.


We can use the parameterized approach developed in preceding sections to generate the dataset for the mapping from $(t_g,\x_c)$ to the optimal feedback control $C(t_g,\x_c)$. Then, a simple FNN trained by the dataset can be used to represent the optimal feedback control. For the current example, the dataset is generated by setting $N$ and $\Delta t$ in Procedure 1 as $6.0\times10^{5}$ and $0.01$, respectively. Then, an FNN with three hidden layers (each of which contains $30$ neurons) is trained to approximate the mapping $(t_g,\x_c)\mapsto C(t_g,\x_c)$. The training is terminated when the mean-squared error between the predicted values and the real values is less than $10^{-6}$. To examine the real-time performance of the trained FNN, it is tested on the MYC-Y6ULY2 CPU at $528$ MHz for $10$ thousand times, showing that the maximum computational time is $0.28$ ms. This computational time is enough for the GNC system of usual satellites. 

\begin{figure}[!htp]
	\centering
	\includegraphics[width=0.8\textwidth]{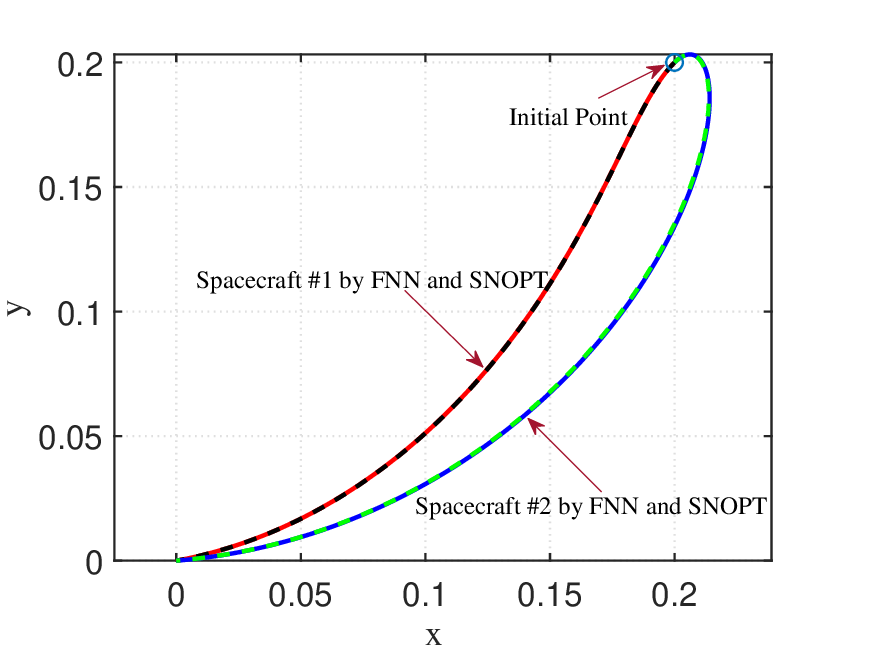}
	\caption{Trajectories for optimal proximity.}
	\label{Fig:optimal_trajectories}
\end{figure}

\begin{figure}[!htp]
	\centering
	\includegraphics[width=0.8\textwidth]{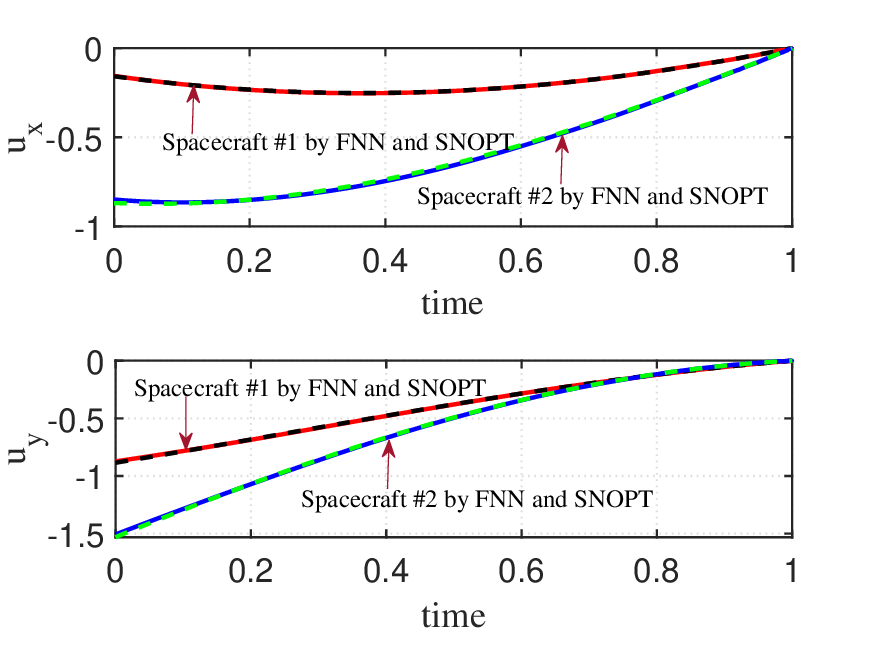}
	\caption{Control profiles for optimal proximity.}
	\label{Fig:optimal_controls}
\end{figure}

\begin{figure}[!htp]
	\centering
	\includegraphics[width=0.8\textwidth]{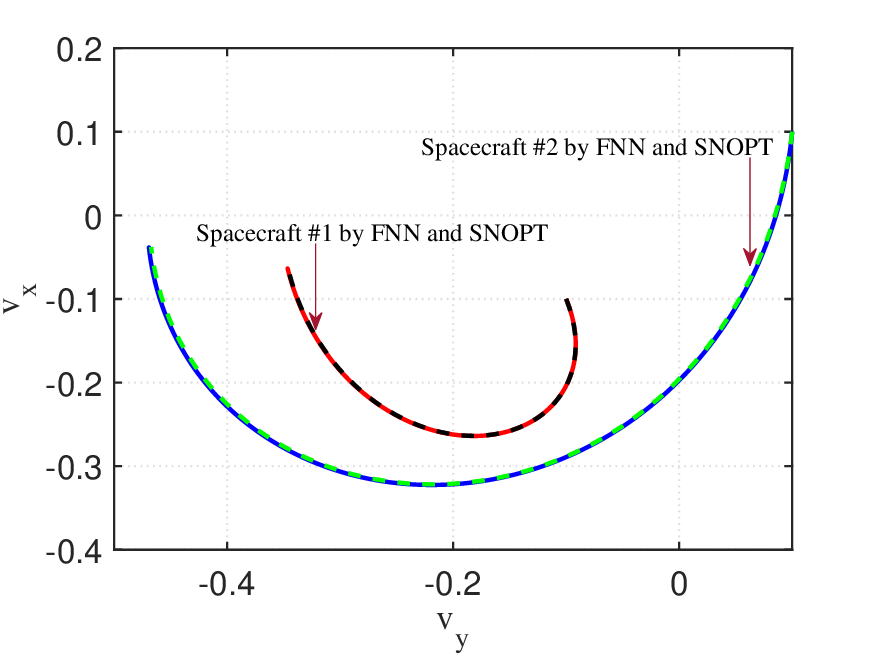}
	\caption{Velocity space profiles for optimal proximity.}
	\label{Fig:optimal_speed}
\end{figure}

We consider two different initial conditions, as shown in Table~ \ref{Tab:initial_condition0}. The trajectories of the two spacecraft are presented in Fig. \ref{Fig:optimal_trajectories}, and their control profiles are demonstrated by Fig. \ref{Fig:optimal_controls}. The velocity space profiles are presented in Fig. \ref{Fig:optimal_speed}. Note that the trajectories by SNOPT (dashed curves) are identical to those by FNN (solid curves). Whereas, it is worth mentioning again here that the optimization methods, like the SNOPT, cannot guarantee to generate optimal trajectories in real time as they suffer the issue of convergence. 

\begin{table}[!htp]
	\centering
	\caption{The initial states and initial time-to-go for each spacecraft.}\label{Tab:initial_condition0}
	\begin{tabular}{cccccc}
		\toprule
		States & $x_0$  & $y_0$ & $v_{x0}$  & $v_{y0}$ & $t_{g0}$\\
		\hline
		Spacecraft $\# 1$ & $0.2$  & $0.2$   & $-0.1$  &  $-0.1$   & 1  \\
		Spacecraft $\# 2$ & $0.2$  & $0.2$   & $0.1$  &  $0.1$   & 1   \\
		\bottomrule
	\end{tabular}
\end{table}


\begin{figure}[!htp]
	\centering
	\includegraphics[width=0.8\textwidth]{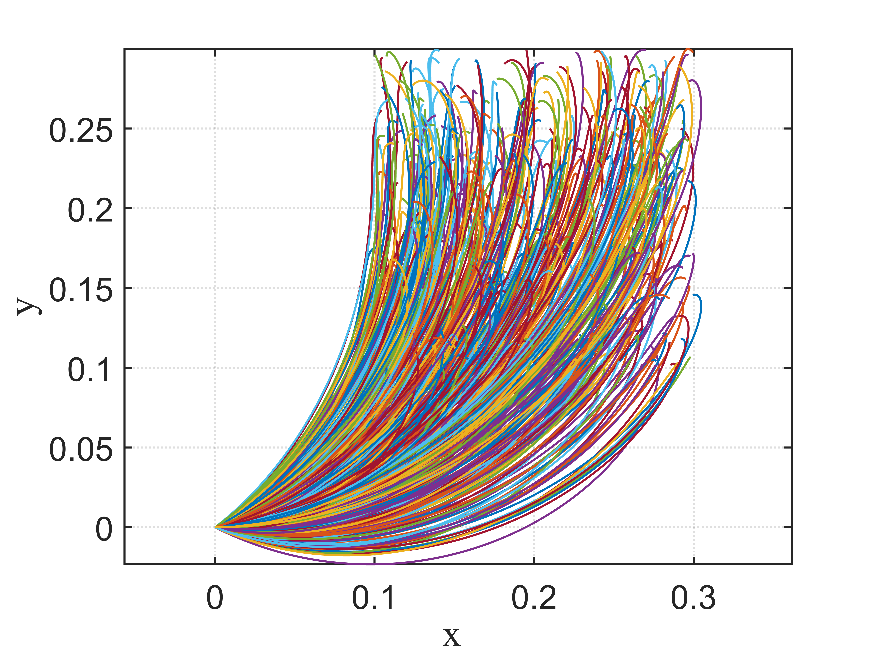}
	\caption{Trajectories of Monte Carlo test for optimal proximity.}
	\label{Fig:monte_trajectories}
\end{figure}

\begin{figure}[htbp]
	\centering
	\subfigure[Errors of terminal tangential displacemant.]{
		\includegraphics[width=0.8\textwidth]{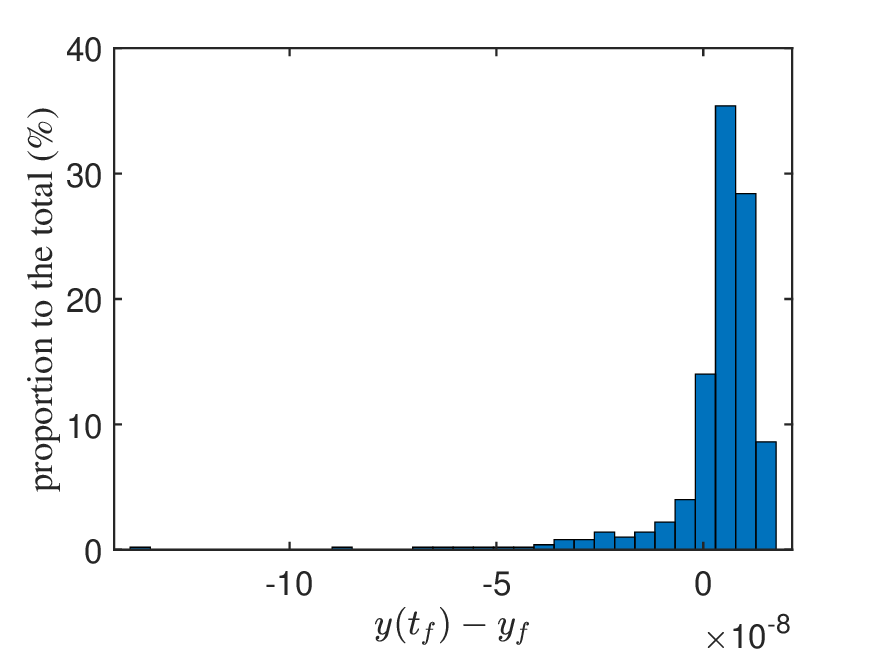}
		\label{Fig:monte_yerror1}
	}
	\quad
	\subfigure[Errors between initial time-to-go and final time.]{
		\includegraphics[width=0.8\textwidth]{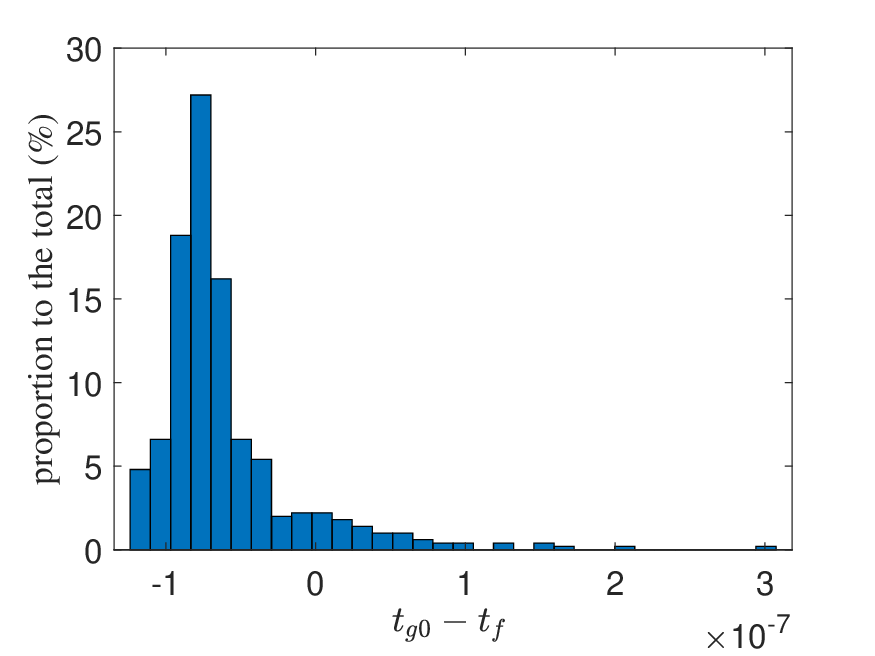}
	}H
	\caption{Histograms for errors of Monte Carlo test of optimal proximity.}
	\label{Fig:monte_error1}
\end{figure}

The Monte Carlo test is also implemented to check the performance of FNN in optimal proximity. $500$ simulations are conducted by randomly sampling the initial conditions $x_0$, $y_0$, $v_{x0}$, and $v_{y0}$ in $[0.1,0.3]$, $[0.1,0.3]$, $[-0.1,0.1]$, $[-0.1,0.1]$, respectively. The trajectories generated in this test are presented in Fig.~\ref{Fig:monte_trajectories}, and the terminal errors are demonstrated by histograms in Fig.~\ref{Fig:monte_error1}.  It is clearly seen from Fig.~\ref{Fig:monte_yerror1} that the maximum error for terminal tangential displacement is less than $1.5\times 10^{-7}$. We consider, for example, that the radius of the target's circular orbit is $10,000$ km. Then, the distribution ranges of $x_0$ and $y_0$ are up to $[1000, 3000]$ km, and $[1000, 3000]$ km, and the distribution ranges of $v_{x0}$ and $v_{y0}$ are up to $[-631.35, 631.35 ]$ m/s and  $[-631.35,631.35]$ m/s, respectively. With such large distributions on the initial conditions, the maximum error for terminal tangential displacement, by converting to the unit of meter, is less than $1.5$ m. 




\section{Conclusions}\label{SE:Conclusions}
The problem of generating optimal flight trajectories for MECPs in real time via ANN was studied in the paper. In order to make sure that an ANN is trained to generate optimal flight trajectories, one usually uses indirect or direct methods to construct the dataset of trajectories for training. However, the trajectories constructed by indirect and direct methods cannot be guaranteed to be at least locally optimal. Instead of using indirect and direct methods, this paper proposed a parametrization approach, which embedding not only necessary conditions but also sufficient conditions for optimality into a set of parameterized differential equations. Furthermore, by establishing the boundary conditions for the parameterized differential equations,  it is enough to solve some initial value problems in order to construct the dataset for the mapping from flight state to optimal feedback control command, as shown by {\it Procedure \ref{algo1}}. Two typical MECPs in aerospace engineering were presented to illustrate how to use the parametrization approach in Procedure 1 to construct the dataset for training FNNs. Numerical simulations first shew that the trained FNNs can generate optimal trajectories in milliseconds, which usually cannot be achieved by indirect and direct methods. In addition, because each trajectory in the dataset constructed by Procedure 1 is at least locally optimal, the trained FNNs sometimes generate trajectories with smaller performance indices, as shown by the example of optimal gliding of flight vehicles.  

\section*{Acknowledgment}

This research was supported by the National Natural Science Foundation of China with grant No. 62088101.


\bibliographystyle{plain}        
 \bibliography{IEEEabrv,Bibliography}                            

\end{document}